\documentclass{amsart}

\usepackage{amssymb}
\usepackage{upref}
\usepackage[arrow,matrix,curve,cmtip,ps]{xy}

\newtheorem{thm}{Theorem}[section]
\newtheorem{cor}[thm]{Corollary}
\newtheorem{lem}[thm]{Lemma}
\newtheorem{prop}[thm]{Proposition}

\newcommand{\secref}[1]{Section~\textup{\ref{#1}}}
\newcommand{\thmref}[1]{Theorem~\textup{\ref{#1}}}
\newcommand{\corref}[1]{Corollary~\textup{\ref{#1}}}
\newcommand{\lemref}[1]{Lemma~\textup{\ref{#1}}}
\newcommand{\propref}[1]{Proposition~\textup{\ref{#1}}}

\theoremstyle{definition}
\newtheorem{defn}[thm]{Definition}
\newtheorem*{warn*}{Warning}

\newcommand{\midtext}[1]{\quad\text{#1}\quad}

\numberwithin{equation}{section}

\newcommand{\mor}{\operatorname{Mor}}
\newcommand{\infl}{\operatorname{Inf}}
\newcommand{\ind}{\operatorname{Ind}}
\newcommand{\ad}{\operatorname{Ad}}
\newcommand{\yg}{Y_{G/G}^G}
\newcommand{\yn}{Y_{G/N}^G}
\newcommand{\zg}{Z_{G/G}^G}
\newcommand{\zn}{Z_{G/N}^G}
\newcommand{\xn}{X_N^G}
\newcommand{\xg}{X_{e}^G}

\renewcommand{\bar}{\overline}
\newcommand{\CC}{\mathcal C}
\newcommand{\A}{\mathcal A}
\newcommand{\K}{\mathcal K}
\renewcommand{\d}{\delta}
\newcommand{\e}{\epsilon}
\renewcommand{\a}{\alpha}
\renewcommand{\b}{\beta}
\newcommand{\g}{\gamma}
\renewcommand{\l}{\lambda}
\renewcommand{\r}{\rho}
\newcommand{\p}{\phi}
\newcommand{\z}{\zeta}
\renewcommand{\P}{\Phi}
\newcommand{\D}{\Delta}
\newcommand{\T}{\Theta}
\renewcommand{\S}{\Sigma}

\newcommand{\id}{\text{id}}
\renewcommand{\:}{\colon}
\newcommand{\<}{\langle}
\renewcommand{\>}{\rangle}
\newcommand{\inv}{^{-1}}
\renewcommand{\hat}{\widehat}
\newcommand{\hathatd}{\hat{\hat\d}}
\newcommand{\hathate}{\hat{\hat\e}}
\newcommand{\hathatz}{\hat{\hat\z}}
\newcommand{\hathata}{\hat{\hat\a}}
\renewcommand{\tilde}{\widetilde}
\newcommand{\dec}{^{\text{dec}}}
\newcommand{\cotimes}{\mathrel{\sharp}}

\title{Mansfield's Imprimitivity Theorem for Full Crossed
Products}

\author{S. Kaliszewski}
\address{Department of Mathematics and Statistics\\
Arizona State University\\
Tempe, Arizona 85287}
\email{kaliszewski@asu.edu}

\author{John Quigg}
\email{quigg@math.asu.edu}



\begin{document}

\begin{abstract}
For any maximal coaction $(A,G,\d)$ and any closed normal
subgroup $N$ of $G$, there exists an
imprimitivity bimodule $Y_{G/N}^G(A)$
between the full crossed product
$A\times_\d G\times_{\hat\d|}N$ and $A\times_{\d|}G/N$,
together with an
$\infl\hat{\hat\d|}-\d\dec$ compatible coaction $\d_Y$ of $G$.
The assignment $(A,\d)\mapsto (Y_{G/N}^G(A),\d_Y)$
implements a natural equivalence between the crossed-product
functors ``${}\times G\times N$'' and ``${}\times G/N$'', in
the category whose objects are maximal coactions of $G$ and
whose morphisms are isomorphism classes of 
right-Hilbert bimodule coactions of $G$.
\end{abstract}

\maketitle

\section{Introduction}

For an action $\alpha$ of a locally compact group $G$ on a
$C^*$-algebra $B$ and a closed normal subgroup $N$ of $G$,
Green's imprimitivity theorem (as formulated, for example,
in \cite[Theorem~B.5]{enchilada}) says that the iterated crossed
product $B\times_\a G\times_{\hat\a|}G/N$ by the restriction
of the dual coaction $\hat\a$ is Morita-Rieffel equivalent
to $B\times_{\a|}N$.   Thus, a representation of $B\times_\a G$
is induced from a representation of $B\times_{\a|}N$ if
and only if it can be combined with a representation of
$C_0(G/N)$ to make a covariant representation of the system
$(B\times_\a G, G/N, \hat\a|)$.
This represents the culmination of 
roughly 50 years of imprimitivity theorems, dating back to
Mackey \cite{mackey}.  

Green's theorem starts with an action, and uses a dual coaction; for
symmetry one desires a ``dual'' version which starts with a coaction and
uses a dual action.
In 1991 Mansfield \cite{man} gave the first such result:
if $\delta$
is a \emph{reduced} coaction of $G$ on a $C^*$-algebra
$A$, and if the normal subgroup
$N$ is \emph{amenable}, then $A\times_\d G\times_{\hat\d|} N$ is
Morita-Rieffel equivalent to $A\times_{\d|} G/N$. This theorem suffers
from two significant drawbacks: 
the coaction is reduced, meaning it uses the
reduced algebra $C^*_r(G)$ and the theory is decidedly spatial; and
the subgroup $N$ is required to be amenable, in order to
have a surjection $C^*_r(G)\to C^*_r(G/N)$ which vouchsafes
the restricted coaction $\d|$.
In \cite{kq:imprimitivity} we removed these
drawbacks by using \emph{full} coactions, 
meaning the full group algebra $C^*(G)$ is used and the theory is
non-spatial; consequently we no longer needed $N$ to be amenable.
However, the price was that we then had to
require our coactions to be \emph{normal} (an unfortunate
term meaning essentially that crossed product duality holds
for reduced crossed products) 
and use the \emph{reduced} crossed by the dual action.

This is in some sense a ridiculous situation; the whole point of crossed
products is to encode \emph{all} the covariant representations,
not just those that factor through the regular representation.
What is of course desired is to use \emph{full} crossed products
throughout. 
The first result in this vein was given in 1998 by 
Echterhoff, Raeburn, and
the first author \cite{ekr}: if $\d$ is a \emph{dual} coaction,
then $A\times_\d G\times_{\d|}
N$ is Morita-Rieffel equivalent to $A\times_{\d|} G/N$, 
where all crossed products are now full. 
Relaxing the restriction that $\d$ be
dual was partially accomplished by Echterhoff and the second
author \cite{eq:full}, but only for \emph{discrete} groups. The missing
ingredient was a crossed product duality for full crossed products. With
Echterhoff, we recently introduced \emph{maximal} coactions \cite{ekq},
which are \emph{by definition} those for which full crossed product
duality holds; 
dual coactions and discrete coactions are specific examples.
The term maximal was chosen because 
every for every coaction $(A,G,\d)$ there is a maximal coaction
$\d^m$ on an extension $A^m$ of $A$ such that
$A^m\times_{\d^m}G\cong A\times_\d G$.
By contrast, 
there is also a normal coaction $\d^n$ on a quotient $A^n$ of $A$
with the same crossed product \cite{qui:fullreduced}
(and in retrospect perhaps a better term for normal coactions
might be \emph{minimal}).

In this paper we prove Mansfield's imprimitivity theorem for full
crossed products (\thmref{full mansfield}):
for every maximal coaction $(A,G,\d)$ and any closed normal
subgroup $N$ of $G$,  the full crossed product
$A\times_\d G\times_{\hat\d|}N$ is Morita-Rieffel equivalent to
$A\times_{\d|} G/N$.  We also prove that this equivalence
is equivariant for the appropriate coactions of $G$.
As mentioned above, an important ingredient in the proof
is a full crossed-product duality theorem for maximal coactions,
which allows us to switch from an arbitrary maximal coaction
to a dual coaction.  The bimodule involved, which we call
the \emph{Katayama bimodule}, is also of interest for non-maximal
coactions, and we investigate it in detail in \secref{kat-sec}.

We also fit our version of Mansfield's
imprimitivity theorem into a categorical perspective that
we, along with Echterhoff and Raeburn \cite{taco, enchilada}, have been
pushing: there is a category whose objects are $C^*$-algebras and whose
morphisms are isomorphism classes of right-Hilbert bimodules (see
\secref{prelims}
for a more detailed summary), and for a fixed group $G$
there are equivariant
versions of this category which incorporate actions and coactions
of $G$. Crossed products are functors among
these categories, and imprimitivity theorems express natural
equivalences between them.
The naturality theorem for full Mansfield imprimitivity
is \thmref{natural mansfield}.
Our proof exploits a curious relationship (\corref{Y-fact})
between Green's imprimitivity bimodule and the Katayama bimodule; 
naturality of those bimodules is proved in Sections~\ref{green-sec}
and~\ref{kat-sec}, respectively.

In \secref{sec:preserve} we prove that maximality is preserved
under various operations: restriction to quotients, inflation
from subgroups, decomposition, and restriction to invariant ideals.  
Besides being of obvious interest
and technical utility, 
the results in this section allow us to 
say that the functors in the
naturality-of-Mansfield theorem
actually live completely within
the category of maximal coactions.

\section{Preliminaries}
\label{prelims}

We adopt the conventions of \cite{enchilada} for almost everything.
However, in some aspects our conventions are slightly different from other
sources, so to avoid any confusion we summarize the possibly unfamiliar
details here.

\subsection{Right-Hilbert bimodules}

Let $A$ and $B$ be $C^*$-algebras.
A \emph{right-Hilbert $A-B$ bimodule},
sometimes denoted ${}_AX_B$,
is a (right) Hilbert $B$-module
$X$ with a compatible, nondegenerate, left $A$-action. We do not
require the $B$-valued inner product to be full, that is, the closed
ideal
\[
B_X:=\bar{\<X,X\>_B}
\]
of $B$ is allowed to be proper.  Compatibility of the left 
$A$-action means that $X$ is an $A-B$ bimodule in the usual algebraic
sense, that is, 
$(a\cdot x)\cdot b=a\cdot (x\cdot b)$ 
for all $a\in A$, $x\in X$, and $b\in B$. 
Nondegeneracy can be interpreted in two equivalent ways: the closed
span $\bar{AX}$ coincides with $X$, or the associated homomorphism
$\kappa_A\:A\to M(\K(X))$ is nondegenerate, where $\K(X)$ is
the imprimitivity algebra of the Hilbert $B$-module $X$.

A right-Hilbert bimodule ${}_AX_B$ is \emph{full} if it is full
as a Hilbert $B$-module, that is, if $B_X=B$.  By restricting the
right-hand coefficient algebra, every right Hilbert $A-B$ bimodule $X$
may be regarded as a full right-Hilbert $A-B_X$ bimodule.  Also, we
may regard $X$ as a right-Hilbert $\K(X)-B$ bimodule.

If $\p\:A\to M(B)$ is a homomorphism, the closed span $\bar{\p(A)B}$ is
called a \emph{standard} right-Hilbert $A-B$ bimodule.  If $\p$ is
nondegenerate, that is, if $\bar{\p(A)B}=B$, then $B$ is called a
\emph{nondegenerate standard} right-Hilbert bimodule.  If ${}_AX_B$ is
any right-Hilbert bimodule, then $\kappa_A$ makes $\K(X)$ into 
a nondegenerate standard $A-\K(X)$ bimodule.

A \emph{partial $A-B$ imprimitivity bimodule} is a right-Hilbert
bimodule ${}_AX_B$ which is also a left Hilbert $A$-module such that
the inner products satisfy the compatibility condition
\[
{}_A\<x,y\>z=x\<y,z\>_B
\]
for all $x,y,z\in X$.
Of course, if both inner products are full then $X$ is an imprimitivity
bimodule in the usual sense.  A \emph{right-partial imprimitivity
bimodule} is a partial imprimitivity bimodule whose left inner product
is full (and similarly for left-partial).  Every right-Hilbert $A-B$
bimodule may be regarded as a $\K(X)-B_X$ imprimitivity bimodule, or
alternatively as a right-partial $\K(X)-B$ imprimitivity bimodule.

Every partial $A-B$ imprimitivity bimodule $X$ has a \emph{linking
algebra} $L(X)$, which is a $C^*$-algebra
characterized by the existence of
complementary projections $p,q\in M(L(X))$ such that
\[
A=pL(X)p,
\quad
B=qL(X)q,
\midtext{and}
X=pL(X)q.
\]
Of course, the projections $p$ and $q$ are full if and only if $X$ is
actually an $A-B$ imprimitivity bimodule.

The \emph{multiplier bimodule} of a right-Hilbert $A-B$ bimodule $X$
is the space $M(X)$ of maps $T\colon B\to X$ which are 
adjointable in the sense that there exists $T^*\colon X\to B$ with
\[
\<Tb,x\>_B = b^*T^*x   
\]
for all $x\in X$ and $b\in B$.  
(Such maps are automatically bounded and $B$-linear.)
$M(X)$ is a right-Hilbert
$M(A)-M(B)$ bimodule in a natural way.  

If ${}_AX_B$ and ${}_CY_D$ are right-Hilbert bimodules and
$\p\colon A\to M(C)$ and $\psi\colon B\to M(D)$ are homomorphisms,
a \emph{$\p-\psi$ compatible right-Hilbert bimodule homomorphism}
is a linear map $\Phi\colon X\to M(Y)$ such that
\[
\Phi(a\cdot x)=\phi(a)\cdot\Phi(x),\quad
\Phi(x\cdot b)=\Phi(x)\cdot\psi(b),\quad
\text{and}\quad
\<\Phi(x),\Phi(z)\>_{M(D)} = \psi(\<x,z\>_B)
\]
for all $a\in A$, $x,z\in X$, and $b\in B$. 
$\Phi$~is \emph{nondegenerate} if $\bar{\Phi(X)D}=Y$ and
the \emph{coefficient maps}
$\p$ and $\psi$ are nondegenerate. 

An \emph{isomorphism} of right-Hilbert $A-B$ bimodules $X$ and $Y$ is 
linear bijection $\Phi\colon X\to Y$ which is an
$\id_A-\id_B$ compatible right-Hilbert bimodule homomorphism.
For the homomorphism property,
it suffices to check the first and third 
conditions above, since
then $\P$ automatically preserves the right
$B$-action.

\subsection{The basic category}

Our basic category $\CC$ has $C^*$-algebras as objects, and the
morphisms from an object $A$ to an object $B$ are the isomorphism
classes of right-Hilbert $A-B$ bimodules.
(We ignore the nuisance of $\mor(A,B)$ not being a set---this never 
causes problems in practice; more precisely, in any given situation of
interest the choice of morphisms can be suitably restricted.)

Given right-Hilbert bimodules ${}_AX_B$ and ${}_BY_C$, the balanced
tensor product $X\otimes_B Y$ is a right Hilbert $A-C$ bimodule.  
The
\emph{composition} $[Y]\circ[X]$ of the morphisms $[X]$ and $[Y]$ 
in the category $\CC$ is the  
isomorphism class $[X\otimes_B Y]$.
We often express the fact that ${}_AZ_C\cong {}_AX\otimes_B Y_C$
by declaring that the diagram
\[
\xymatrix{
{A}
\ar[rr]^{Z}
\ar[dr]_{X}
&&{C}
\\
&{B}
\ar[ur]_{Y}
}
\]
commutes.
The \emph{identity morphism} on an object $A$ is the isomorphism class
of the standard right-Hilbert bimodule ${}_AA_A$.
The \emph{equivalences} (that is, the invertible
morphisms) in $\CC$ are exactly the isomorphism classes of
imprimitivity bimodules.

Every morphism in $\CC$ can be factored as a composition of a
nondegenerate standard morphism and a right-partial equivalence:
namely,
\[
{}_AX_B\cong {}_A\K(X)\otimes_{\K(X)} X_B,
\]
so $[{}_AX_B]=[{}_{\K(X)} X_B]\circ[{}_A\K(X)_{\K(X)}]$.

\subsection{Equivariant categories}

Given a right-Hilbert bimodule ${}_AX_B$ and actions $(A,\a)$ and 
$(B,\b)$ of a
locally compact group $G$,
an \emph{$\a-\b$ compatible action} on $X$ is a point-norm
continuous map $\g$ of $G$ into the linear bijections on $X$ such
that
\[
\g_s(ax)=\a_s(a)\g_s(x)
\quad\text{and}\quad
\<\g_s(x),\g_s(y)\>_B
=\b_s\bigl(\<x,y\>_B\bigr)
\]
for all $s\in G$, $a\in A$, and $x,y\in X$. 
(The property $\g_s(xb)=\g_s(x)\b_s(b)$
for all $s\in G$, $x\in X$, and $b\in B$
is then automatic.)  

Given another action $(C,\e)$ and a
$\b-\e$ compatible right-Hilbert bimodule action $(Y,\r)$, the $\a-\e$
compatible \emph{balanced tensor product action} $\g\otimes_B\r$ on
$X\otimes_BY$ is given by
\[
(\g\otimes_B\r)_s(x\otimes y)=\g_s(x)\otimes\r_s(y).
\]

If $(X,\g)$ and $(Y,\r)$ are $\a-\b$ compatible actions, a
right-Hilbert bimodule isomorphism $\P\:X\to Y$ is \emph{$\g-\r$
equivariant} if it intertwines $\g$ and $\r$, that is, if
$\P\circ\g_s=\r_s\circ\P$ for all $s\in G$.  

The category $\A(G)$ of actions of $G$ has $C^*$-algebra
actions as
objects and isomorphism classes of right-Hilbert bimodule actions as
morphisms.  The composition is given by the balanced tensor product:
\[
[(Y,\r)]\circ[(X,\g)]=[(X\otimes_BY,\g\otimes_B\r)],
\]
and the identity morphism on an object $(A,\a)$ is the isomorphism
class of the standard bimodule action ${}_{(A,\a)}(A,\a)_{(A,\a)}$.

Turning to the dual category, we must first issue a warning:
\begin{warn*}
We assume without further
comment that all $C^*$-coactions $(A,G,\d)$ satisfy
the nondegeneracy condition
\begin{equation}\label{nondegen}
\bar{\d(A)(1\otimes C^*(G))}=A\otimes C^*(G).
\end{equation}
\end{warn*}
It is still an open problem whether the
above nondegeneracy condition is automatic, and if it should turn out
that non-nondegenerate coactions exist, they will be of no interest in
duality theory.  
Thus, we feel that the nondegeneracy condition should
be included in the definition of coaction.

Given a right-Hilbert bimodule ${}_AX_B$
and coactions $(A,\d)$ and $(B,\e)$ of $G$, a \emph{$\d-\e$ compatible
coaction} on $X$ is a nondegenerate $\d-\e$ compatible
right-Hilbert bimodule homomorphism
$\z\:X\to M(X\otimes C^*(G))$ which satisfies the nondegeneracy condition
\[
\bar{(1\otimes C^*(G))\z(X)}=X\otimes C^*(G)
\]
and the coaction identity
\[
(\z\otimes\id)\circ\z=(\id\otimes\d_G)\circ\z.
\]
Thus we have built nondegeneracy of $\z$, both as a bimodule
homomorphism \emph{and as a coaction}, into the definition.  The
companion nondegeneracy condition
\[
\bar{\z(X)(1\otimes C^*(G)}=X\otimes C^*(G)
\]
is then automatic.

Given another coaction $(C,\vartheta)$ and a $\e-\vartheta$ compatible
right-Hilbert bimodule coaction $(Y,\eta)$, the $\d-\vartheta$
compatible \emph{balanced tensor product coaction} $\z\cotimes_B\eta$
on $X\otimes_BY$ is given by
\[
\z\cotimes_B\eta=\T\circ(\z\otimes_B\eta),
\]
where $\T\:(X\otimes C^*(G))\otimes_{B\otimes C^*(G)}(Y\otimes
C^*(G))\to (X\otimes_BY)\otimes C^*(G)$ is the isomorphism given by
\[
\T\bigl((x\otimes c)\otimes(y\otimes d)\bigr)
=(x\otimes y)\otimes cd.
\]

If $(X,\z)$ and $(Y,\eta)$ are $\d-\e$ compatible coactions, a
right-Hilbert bimodule isomorphism $\P\:X\to Y$ is \emph{$\z-\eta$
equivariant} if
\[
\eta\circ\P=(\P\otimes\id)\circ\z.
\]

The category $\CC(G)$ of coactions of $G$ has
$C^*$-algebra coactions as objects and isomorphism classes of 
right-Hilbert
bimodule coactions as morphisms. The composition is given by
the balanced
tensor product:
\[
[(Y,\eta)]\circ[(X,\z)]=[(X\otimes_BY,\z\cotimes_B\eta)],
\]
and the identity morphism on an object $(A,\d)$ is the isomorphism
class of the standard bimodule coaction ${}_{(A,\d)}(A,\d)_{(A,\d)}$.

The familiar operations of forming crossed products, inflating,
restricting, and forming decomposition actions or coactions are functors
among the equivariant categories. 

\section{Naturality of Green's Bimodule}
\label{green-sec}

Recall that for any action $(B,G,\a)$ and any closed normal
subgroup $N$ of $G$, Green's imprimitivity theorem gives a 
$B\times_\a G\times_{\hat\a|} G/N-B\times_{\a|} N$
imprimitivity bimodule $X_N^G(B)$
which is a completion of $C_c(G,B)$
(\cite{gre:local}; see also \cite[Appendix~B]{enchilada}).
There is also a
$\hat\a\dec-\infl\hat{\a|}$
compatible coaction $\d_{X(B)}$ on $X_N^G(B)$ which
maps $C_c(G,B)$ into $C_c(G,M^\beta(B\otimes C^*(G)))
\subset M(X_N^G(B)\otimes C^*(G))$  by the formula
\[
\d_{X(B)}(x)(t) = x(t)\otimes t
\]
\cite[Section~4.4]{enchilada}. 
(The superscript $\b$ means we require the 
functions to be continuous for the
strict topology on the multipliers.)

In this section we prove that Green's bimodule is a 
coaction-equivariant natural
equivalence between certain crossed product functors. 
Without the equivariance, this was proven in \cite{taco}.
With equivariance, it was  proven in 
\cite[Theorem 4.20]{enchilada}, but in the
context of reduced crossed products, and
there the construction was
deduced from a more general one involving induced actions. 
Since we rely heavily on naturality of Green's bimodule
for some of our main results,
we feel it is best to give the details here. 
It also provides a good illustration of the 
techniques which are typical in proving naturality theorems.

\begin{thm}[Naturality of Green]
\label{natural green}
Let ${}_{(B,\a)}(X,\g)_{(C,\b)}$ be a right-Hilbert bimodule
action of a locally compact group $G$, and let $N$ be a 
closed normal subgroup of $G$.  Then the diagram
\[
\xymatrix
@C+30pt
{
{(B\times_\a G\times_{\hat\a|} G/N,\hat\a\dec)}
\ar[r]^-{(X_N^G(B),\d_{X(B)}}
\ar[d]_{(X\times G\times G/N,\hat\g\dec)}
&{(B\times_{\a|} N,\infl\hat{\a|})}
\ar[d]^{(X\times N,\infl\hat{\g|})}
\\
{(C\times_\b G\times_{\hat\b|} G/N,\hat\b\dec)}
\ar[r]_-{(X_N^G(C),\d_{X(C)})}
&{(C\times_{\b|} N,\infl\hat{\b|})}
}
\]
commutes.
\end{thm}

\begin{proof}
We follow the strategy employed in both \cite{taco} 
and \cite{enchilada}:
any such $(X,\g)$ can be factored as a 
balanced tensor product of a
nondegenerate standard bimodule action and a right-partial
imprimitivity bimodule action.  Once we know the commutativity for
each of these two special cases, the result follows by chaining
together the corresponding diagrams, using 
functoriality of the crossed product.

So, first assume $(X,\g)$ is of the form $(C,\b)$, where we have a
nondegenerate $\a-\b$ equivariant homomorphism $\p\:B\to M(C)$.%
\footnote{For this part of the proof, by \cite[Lemma~4.10]{enchilada}
it suffices to find an equivariant imprimitivity bimodule
homomorphism of $\xn(B)$ onto $\xn(C)$ whose coefficient maps
are $\p\times G\times G/N$ and $\p\times N$.  
Instead, we will take advantage of the bimodule isomorphism
already provided by \cite{taco}.}
The required diagram is the outer rectangle of the diagram
\[
\xymatrix
@C+30pt
{
{(B\times_\a G\times_{\hat\a|} G/N,\hat\a\dec)}
\ar[r]^-{(X_N^G(B),\d_{X(B)})}
\ar[d]_{\p\times G\times G/N}
\ar[dr]|{(X_N^G(C),\d_{X(C)})}
&{(B\times_{\a|} N,\infl\hat{\a|})}
\ar[d]^{\p\times N}
\\
{(C\times_\b G\times_{\hat\b|} G/N,\hat\b\dec)}
\ar[r]_-{(X_N^G(C),\d_{X(C)})}
&{(C\times_{\b|} N,\infl\hat{\b|})},
}
\]
where the bimodule on the diagonal is defined by the lower left
triangle; the coaction on the diagonal is compatible with
the appropriate coactions precisely because 
$\p\times G\times G/N$
is $\hat\a\dec-\hat\b\dec$ equivariant.

For the upper triangle, we use the right-Hilbert
bimodule isomorphism
$\Psi$ of $X_N^G(B)\otimes_{B\times N}(C\times_{\b|} N)$
onto $X_N^G(C)$
provided by \cite{taco}. 
For an elementary tensor $x\otimes g\in
C_c(G,B)\odot C_c(N,C)$, 
$\Psi(x\otimes g)\in C_c(G,C)$ is defined by 
\[
\Psi(x\otimes g)(t)=\int_N\p(x(tn))\b_{tn}(g(n\inv))\D(n)^{-1/2}\,dn.
\]
To see that $\Psi$ is 
$\d_{X(B)}\cotimes_{B\times N} \infl\hat{\b|}-\d_{X(C)}$
equivariant, compute:
\begin{align}
&(\Psi\otimes\id)\circ(\d_{X(B)}\cotimes \infl\hat{\b|})(x\otimes g)(t)
\notag
\\&\quad=(\Psi\otimes\id)\circ\T\bigl(
\d_{X(B)}(x)\otimes \infl\hat{\b|}(g)\bigr)(t)
\notag
\\&\quad=\int(\p\otimes\id)\bigl(\d_{X(B)}(x)(tn)\bigr)
(\b_{tn}\otimes\id)\bigl(\infl\hat{\b|}(g)(n\inv)\bigr)
\D(n)^{-1/2}\,dn
\tag{*}
\\&\quad=\int(\p\otimes\id)(x(tn)\otimes tn)
(\b_{tn}\otimes\id)(g(n\inv)\otimes n\inv)
\D(n)^{-1/2}\,dn
\notag
\\&\quad=\int\bigl(\p(x(tn))\otimes tn\bigr)
\bigl(\b_{tn}(g(n\inv))\otimes n\inv\bigr)
\D(n)^{-1/2}\,dn
\notag
\\&\quad=\int\bigl(\p(x(tn))\b_{tn}(g(n\inv))
\otimes t\bigr)\D(n)^{-1/2}\,dn
\notag
\\&\quad=\int\p(x(tn))\b_{tn}(g(n\inv))\D(n)^{-1/2}\,dn
\notag
\otimes t
\\&\quad=\Psi(x\otimes g)(t)\otimes t
\notag
\\&\quad=\d_{X(C)}\circ\Psi(x\otimes g)(t).
\notag
\end{align}
The equality (*) is verified by temporarily
replacing $\d_{X(B)}(x)$ and $\infl\hat{\b|}(g)$ by elementary tensors
$x'\otimes c\in C_c(G,B)\odot C^*(G)$ and $g'\otimes d\in C_c(N,C)\odot
C^*(G)$, computing that 
\begin{align*}
&(\Psi\otimes\id)\circ\T\bigl(
(x'\otimes c)\otimes(g'\otimes d)\bigr)(t)
\\&\quad=(\Psi\otimes\id)\bigl((x'\otimes g')\otimes cd\bigr)(t)
\\&\quad=\bigl(\Psi(x'\otimes g')\otimes cd\bigr)(t)
\\&\quad=\Psi(x'\otimes g')(t)\otimes cd
\\&\quad=\int\p(x'(tn))\b_{tn}(g'(n\inv))\D(n)^{-1/2}\,dn\otimes cd
\\&\quad=\int\bigl(\p(x'(tn))\b_{tn}(g'(n\inv))\otimes cd\bigr)
\D(n)^{-1/2}\,dn
\\&\quad=\int\bigl(\p(x'(tn))\otimes c\bigr)
\bigl(\b_{tn}(g'(n\inv))\otimes d\bigr)
\D(n)^{-1/2}\,dn
\\&\quad=\int(\p\otimes\id)\bigl(x'(tn)\otimes c\bigr)
(\b_{tn}\otimes\id)\bigl(g'(n\inv)\otimes d\bigr)
\D(n)^{-1/2}\,dn
\\&\quad=\int(\p\otimes\id)\bigl((x'\otimes c)(tn)\bigr)
(\b_{tn}\otimes\id)\bigl((g'\otimes d)(n\inv)\bigr)
\D(n)^{-1/2}\,dn,
\end{align*}
and then appealing to density and continuity for the
inductive-limit topologies.

The argument for the lower triangle is much easier: 
if 
\[
\Upsilon\colon
C\times_\b G\times_{\hat\b|}G/N\otimes_{C\times G\times G/N}
X_N^G(C) \to X_N^G(C)
\]
is the canonical isomorphism, then for 
$f\in C\times_\b G\times_{\hat\b|}G/N$ 
and $x\in X_N^G(C)$ it is not hard to see that
\[
(\Upsilon\otimes\id)\circ(\hat\b\dec\cotimes_{C\times G\times
G/N} \d_{X(C)})(f\otimes x)
= \hat\b\dec(f)\cdot\d_{X(C)}(x),
\]
while
\[
\d_{X(C)}\circ\Upsilon(f\otimes x)
= \d_{X(C)}(f\cdot x);
\]
so it amounts to the fact that $\d_{X(C)}$ is a right-Hilbert
bimodule homomorphism with left coefficient map $\hat\b\dec$.

Now assume that ${}_{(B,\a)}(X,\g)_{(C,\b)}$ is a
right-partial imprimitivity bimodule action.  Let $L(X)$ be the
linking algebra of ${}_BX_C$.  Then there is an action $\nu$ on
$L(X)$ which compresses on the corners to the given actions on $B$,
$X$, and $C$.  It follows straight from the definitions that the
coaction $\d_{X(L(X))}$ on the associated Green bimodule $X_N^G(L(X))$
compresses on the diagonal corners to $\d_{X(B)}$ and $\d_{X(C)}$. 
Using the identification
\[
L(X)\times_\nu G\times_{\hat\nu|}G/N
=L(X\times_\g G\times_{\hat\g|}G/N),
\]
we can apply the linking algebra techniques of
\cite[Section 4.2]{enchilada} to
conclude that the desired diagram commutes, provided we can show
\[
\<pX_N^G(L(X))q,pX_N^G(L(X))q\>_{C\times N}
\subset \bar{\<X\times N,X\times N\>_{C\times N}}.
\]
But we have
\begin{align*}
\<pC_c(G,L(X))q,pC_c(G,L(X))q\>_{C\times N}
&=\<C_c(G,X),C_c(G,X)\>_{C\times N}
\\&\subset \bar{\<X\times N,X\times N\>_{C\times N}},
\end{align*}
which implies the desired containment since $pC_c(G,L(X))q$ is dense in
$pX_N^G(L(X))q$.
\end{proof}

\section{The Katayama bimodule}
\label{kat-sec}

In this section we put no restrictions on our coaction 
(except that we continue to assume it satisfies the nondegeneracy 
condition of \eqref{nondegen}!), 
and in this generality we introduce
a right-Hilbert bimodule, which we
call the Katayama bimodule, connecting
the double dual coaction to the original; in the next section we
restrict our attention to maximal coactions, and then the Katayama
bimodule will be an imprimitivity bimodule.

Let $(A,G,\d)$ be a coaction. Then
as observed in \cite[Corollary 2.6]{nil:duality}
there is always a \emph{canonical surjection}
\[
\P_A
=(\id\otimes\l)\circ\d\times(1\otimes M)\times(1\otimes\r)
\:A\times_\d G\times_{\hat\d}G\to A\otimes \K(L^2(G)),
\]
where $\l$ and $\r$ are the left and right
regular representations of $G$, and $M$ is
the multiplication representation of $C_0(G)$, on $L^2(G)$.
When there is no danger of confusion we will write 
$\P$ for $\P_A$. Also we will 
usually write $L^2$ for $L^2(G)$, and $\K$ for 
the compact operators on $L^2(G)$.

The external bimodule tensor product
$A\otimes L^2$ is an $A\otimes\K - A$ imprimitivity bimodule,
hence becomes a right-Hilbert 
$A\times_\d G\times_{\hat\d} G - A$ 
bimodule when the left coefficient action is modified via $\P$. 
More precisely:

\begin{defn}
For any coaction $(A,G,\d)$, the \emph{Katayama bimodule}
$K(A)$ is the right-Hilbert $A\times_\d G\times_{\hat\d}G-A$ 
bimodule defined as the Hilbert $A$-module $A\otimes L^2$ with left 
$A\times_\d G\times_{\hat\d}G$-module action given by
\[
f\cdot x=\P(f)x.
\]
\end{defn}

We will now show that there is a $\hathatd-\d$ compatible
coaction $\d_K$ on $K(A)$.
First we recall some notation: $w_G$ is the unitary element of 
$M(C_0(G)\otimes C^*(G))$ defined by the canonical embedding 
$G\hookrightarrow M(C^*(G))$. For any $C^*$-algebras $C$ and $D$, 
$\S\:C\otimes D\to D\otimes C$ is the flip isomorphism. We let 
$\d\otimes_*\id=(\id\otimes\S)\circ(\d\otimes\id)$ denote the 
``external'' tensor product coaction on either 
the bimodule $A\otimes L^2$ or  the $C^*$-algebra $A\otimes\K$. 
Further let
\[
W=1\otimes (M\otimes\id)(w_G^*)\in M(A\otimes\K\otimes C^*(G)),
\]
where the first $M$ is the multiplication representation of $C_0(G)$.

\begin{prop}\label{kat-coact}
For any coaction $(A,G,\d)$, the map
\[
\d_K=W(\d\otimes_*\id_{L^2})
\]
is a $\hathatd-\d$ compatible coaction on the Katayama bimodule
$K(A)$.
\end{prop}

\begin{proof}
The canonical surjection $\P$ 
is $\hathatd-\ad W\circ(\d\otimes_*\id)$ equivariant, hence
$\ad W\circ(\d\otimes_*\id)$ is 
a $\hathatd-\ad W\circ(\d\otimes_*\id)$ compatible 
coaction on the standard right-Hilbert
bimodule ${}_{A\times_\d G\times_{\hat\d}G}(A\otimes\K)_{A\otimes\K}$.
Moreover,
$W$ is a $(\d\otimes_*\id)$-cocycle, so $W(\d\otimes_*\id)$ is 
an $\ad W\circ(\d\otimes_*\id)-(\d\otimes_*\id)$ compatible coaction 
on ${}_{A\otimes\K}A\otimes\K_{A\otimes\K}$.
Also,
$\d\otimes_*\id_{L^2}$ is a $(\d\otimes_*\id_\K)-\d$ compatible
coaction on the imprimitivity bimodule
${}_{A\otimes\K}(A\otimes L^2)_A$.
Combining these, we see that 
the balanced tensor product coaction
\[
\ad W\circ(\d\otimes_*\id_\K)\cotimes_{A\otimes\K}
W(\d\otimes_*\id_\K)\cotimes_{A\otimes\K}(\d\otimes_*\id_{L^2})
\]
on
$(A\times_\d G\times_{\hat\d}G)\otimes_{A\otimes\K}
(A\otimes\K)\otimes_{A\otimes\K}(A\otimes L^2)$
is $\hathatd-\d$ compatible. 

To complete the proof, it suffices to 
show that the canonical right-Hilbert bimodule 
isomorphism
\[
(A\otimes\K)\otimes_{A\otimes\K}
(A\otimes\K)\otimes_{A\otimes\K}(A\otimes L^2)
\to K(A)
\]
takes the above triple tensor product coaction to the map
$\d_K$.
In other words, we want to show that the 
outer rectangle of the diagram
\[
\xymatrix
@C+40pt
{
{(A\times_\d G\times_{\hat\d} G,\hathatd)}
\ar[r]^-{(K(A), \d_K)}
\ar[d]_{\P}
&{(A,\d)}
\\
{(A\otimes\K,\ad W\circ(\d\otimes_*\id))}
\ar[ur]|{(A\otimes L^2,W(\d\otimes_*\id))}
\ar[r]_-{(A\otimes\K,W(\d\otimes_*\id))}
&{(A\otimes\K,\d\otimes_*\id)}
\ar[u]|{(A\otimes L^2,\d\otimes_*\id)}
}
\]
commutes.

Let 
$\Psi\:(A\otimes\K)\otimes_{A\otimes\K}(A\otimes L^2)\to A\otimes L^2$
be the
canonical isomorphism given by $\Psi(a\otimes k\otimes b\otimes
\xi)=(ab\otimes k\xi)$ for $a,b\in A$, $k\in\K$, and $\xi\in L^2$. 
Then we have
\begin{align}
&(\Psi\otimes\id)\circ\bigl(W(\d\otimes_*\id)\cotimes_{A\otimes\K}
(\d\otimes_*\id)\bigr)\bigl((a\otimes k)\otimes (b\otimes\xi)\bigr)
\notag
\\&=(\Psi\otimes\id)\circ\T\bigl(
W(\d\otimes_*\id)(a\otimes k)\otimes(\d\otimes_*\id)(b\otimes\xi)\bigr)
\notag
\\&=W(\d\otimes_*\id)(a\otimes k)(\d\otimes_*\id)(b\otimes\xi)
\tag{*}
\\&=W(\d\otimes_*\id)\bigl((a\otimes k)(b\otimes\xi)\bigr)
\notag
\\&=W(\d\otimes_*\id)\circ\Psi
\notag
\bigl((a\otimes k)\otimes (b\otimes\xi)\bigr)
\end{align}
where the equality (*) is easily verified on elementary tensors: for
$a',b'\in A$, $k'\in \K$, $\xi'\in L^2$, and $c,d\in C^*(G)$ we have
\begin{align*}
(\Psi\otimes\id)\circ\T\bigl(
(a'\otimes k'\otimes c)\otimes(b'\otimes\xi\otimes d)\bigr)
&=(\Psi\otimes\id)\bigl(
(a'\otimes k')\otimes(b'\otimes\xi)\otimes cd\bigr)
\\&=a'b'\otimes k'\xi\otimes cd
\\&=(a'\otimes k'\otimes c)(b'\otimes\xi\otimes d).
\end{align*}
This implies that 

Commutativity of the upper triangle is much simpler:
as in the proof of \thmref{natural green}, 
it amounts to the fact that $W(\d\otimes_*\id)$ is a
right-Hilbert bimodule homomorphism on $A\otimes L^2$ 
with left coefficient map $\ad W\circ(\d\otimes_*\id)$.
\end{proof}

We next show that the Katayama bimodule
is a natural transformation from the functor 
``${}\times G\times G$'' to the identity functor.  

\begin{thm}[Naturality of Katayama]
\label{natural katayama}
Let 
${}_{(A,\d)}(X,\z)_{(B,\e)}$ be a right-Hilbert bimodule coaction
of a locally compact group $G$. 
Then the diagram
\[
\xymatrix
@C+30pt
{
{(A\times_\d G\times_{\hat\d}G,\hathatd)}
\ar[r]^-{(K(A),\d_{K(A)})}
\ar[d]_{(X\times_\z G\times_{\hat\z}G,\hathatz)}
&{(A,\d)}
\ar[d]^{(X,\z)}
\\
{(B\times_\e G\times_{\hat\e}G,\hathate)}
\ar[r]_-{(K(B),\d_{K(B)})}
&{(B,\e)}
}
\]
commutes.
\end{thm}

\begin{proof}
As in the proof of \thmref{natural green}, we use the 
factorization strategy of \cite{taco} and \cite{enchilada}.
So, first assume $(X,\z)$ is of the form $(B,\e)$, where we have a
nondegenerate $\d-\e$ equivariant homomorphism $\p\:A\to M(B)$.  
(Note that we cannot appeal to \cite[Lemma~4.10]{enchilada}
here because the Katayama bimodules are not imprimitivity
bimodules.)
The required diagram is the outer rectangle of the diagram
\[
\xymatrix
@C+30pt
{
{(A\times_\d G\times_{\hat\d}G,\hathatd)}
\ar[r]^-{(K(A),\d_{K(A)})}
\ar[dr]|{(B\otimes L^2,W(\e\otimes_*\id))}
\ar[d]_{\p\times G\times G}
&{(A,\d)}
\ar[d]^{\p}
\\
{(B\times_\e G\times_{\hat\e}G,\hathate)}
\ar[r]_-{(K(B),\d_{K(B)})}
&{(B,\e)},
}
\]
where the bimodule on the diagonal is defined by the
lower triangle;
the coaction on the diagonal is appropriately compatible 
by equivariance of $\phi\times G\times G$.

For the upper triangle, let 
$\Psi\:K(A)\otimes_AB\to B\otimes L^2$ 
be the isomorphism given by
$\Psi(a\otimes\xi\otimes b)=\p(a)b\otimes\xi$.  Then
\begin{align}
&(\Psi\otimes\id)\circ\bigl(
\d_{K(A)}\cotimes_A\e\bigl)(a\otimes\xi\otimes b)
\notag
\\&\quad=(\Psi\otimes\id)\circ\T\bigl(
W(\d\otimes_*\id)(a\otimes\xi)\otimes\e(b)\bigr)
\notag
\\&\quad=W(\id\otimes\S)\bigl(
(\p\otimes\id)(\d(a))\e(b)\otimes\xi\bigr)
\tag{*}
\\&\quad=W(\id\otimes\S)\bigl(
\e(\p(a))\e(b)\otimes\xi\bigr)
\notag
\\&\quad=W(\id\otimes\S)\bigl(
\e(\p(a)b)\otimes\xi\bigr)
\notag
\\&\quad=W(\e\otimes_*\id)(\phi(a)b\otimes\xi)
\notag
\\&\quad=W(\e\otimes_*\id)\circ\Psi(a\otimes\xi\otimes b).
\notag
\end{align}
The equality (*) is verified by first computing, 
for $a'\in A$, $\xi'\in L^2$, $b'\in B$, and $c,d\in C^*(G)$,
that
\begin{align*}
&(\Psi\otimes\id)\circ\T\bigl(
(a'\otimes\xi'\otimes c)\otimes(b'\otimes d)\bigr)
\\&\quad=(\Psi\otimes\id)(a'\otimes\xi'\otimes b'\otimes cd)
\\&\quad=\p(a')b'\otimes\xi'\otimes cd
\\&\quad=(\p(a')\otimes\xi'\otimes c)(b'\otimes 1\otimes d)
\\&\quad=(\p\otimes\id\otimes\id)(a'\otimes\xi'\otimes c)
(\id\otimes\S)(b'\otimes d\otimes 1),
\end{align*}
and then using
\begin{align*}
(\p\otimes\id\otimes\id)\bigl(W(\d\otimes_*\id)(a\otimes\xi)\bigr)
&=W(\p\otimes\id\otimes\id)\bigl(
(\id\otimes\S)(\d(a)\otimes\xi)\bigr)
\\&=W(\id\otimes\S)\bigl((\p\otimes\id)(\d(a))\otimes\xi\bigr)
\end{align*}
to see that
\begin{align*}
&(\Psi\otimes\id)\circ\T\bigl(
W(\d\otimes_*\id)(a\otimes\xi)\otimes\e(b)\bigr)
\\&\quad=(\p\otimes\id\otimes\id)\bigl(
W(\d\otimes_*\id)(a\otimes\xi)\bigr)
(\id\otimes\S)(\e(b)\otimes 1)
\\&\quad=W(\id\otimes\S)\bigl((\p\otimes\id)(\d(a))\otimes\xi\bigr)
(\id\otimes\S)(\e(b)\otimes 1)
\\&\quad=W(\id\otimes\S)\bigl(
((\p\otimes\id)(\d(a))\otimes\xi)(\e(b)\otimes 1)\bigr)
\\&\quad=W(\id\otimes\S)\bigl(
(\p\otimes\id)(\d(a))\e(b)\otimes\xi\bigr).
\end{align*}

Commutativity of the lower triangle amounts to the fact
that $\d_{K(B)}$ is a right-Hilbert bimodule homomorphism
with left coefficient map $\hathate$.

Now assume that ${}_{(A,\d)}(X,\z)_{(B,\e)}$ is a
right-partial imprimitivity bimodule coaction.  
we can apply the linking algebra techniques of 
\cite[Section~4.2]{enchilada} to 
conclude that the desired diagram commutes, provided we can show
\[
\<pK(L(X))q,pK(L(X))q\>_B\subset \bar{\<X,X\>_B}.
\]
But we have
\[
pK(L(X))q
=p(L(X)\otimes L^2)q
=X\otimes L^2,
\]
hence
\[
\<pK(L(X))q,pK(L(X))q\>_B
=\<X\otimes L^2,X\otimes L^2\>_B
=\<X,X\>_B,
\]
so we are done.
\end{proof}


\section{Full Mansfield Imprimitivity}
\label{mans-sec}

In this section, we prove an equivariant version of Mansfield's
Imprimitivity Theorem for full crossed products and maximal
coactions.  Our starting point is 
\cite[Proposition 1.1]{ekr}, which provides the appropriate
imprimitivity bimodule for dual coactions.  Specifically,
given an action $(B,G,\alpha)$ and a closed normal subgroup
$N$ of $G$, there is a
$B\times_\a G\times_{\hat\a}G\times_{\hathata|}N -
B\times_\a G\times_{\hat\a|}G/N$
imprimitivity bimodule $Z_{G/N}^G(B\times_\a G)$.
This bimodule is 
a completion of $C_c(G\times G,B)$, with left action
of $C_c(N\times G\times G, B)$ and right
$C_c(G\times G/N,B)$-valued inner product given by
\begin{align}\label{Z-ops}
f\cdot z(s,t) &= \int_N\int_G
f(n,r,t)\a_r(z(r\inv s,r\inv tn))\,
\Delta_N(n)^{1/2}\, dr\, dn,\\
\notag
\<{z},{w}\>_{B\times G\times G/N}(s,tN) &= \int_N \int_G
\a_{r\inv}(z(r,rtn)^* w(rs,rtn))\, dr\, dn.
\end{align}
This left action is given on generators
$b\in B$, $r\in G$, $h\in C_c(G)\subseteq C_0(G)$,
and $n\in N$ by
\begin{align}\label{Z-gens}
(b\cdot z)(s,t) &= b\,z(s,t) &
(r\cdot z)(s,t) &= \a_r(z(r\inv s,r\inv t))\\\notag
(h\cdot z)(s,t) &= h(t)\, z(s,t)&
(n\cdot z)(s,t) &= z(s,tn)\,\Delta(n)^{1/2}.
\end{align}

The next step is to define an
$\infl\hat{\hathata|}-\hat\a\dec$
compatible coaction $\d_Z$ on $Z_{G/N}^G(B\times_\a G)$. 
A related construction appears in  \cite[Lemma~5.7]{enchilada},
where it was deduced from \cite{er:stab}, and the context there 
is reduced crossed products.  
Thus we feel
it is best to construct the bimodule coaction from scratch here,
and this will take a little bit of work. 

In the following lemma, we
identify $C_c(G\times G,B)\odot C^*(G)$ with
an inductive-limit-dense subspace of 
$C_c(G\times G,B\otimes C^*(G))$, 
and hence with a dense subspace of the bimodule 
$Z_{G/N}^G((B\otimes C^*(G))\times_{\a\otimes\id}G)$, 
by identifying an elementary tensor $z\otimes c$ 
with the function defined by
\[
(z\otimes c)(s,t)=z(s,t)\otimes c.
\]

\begin{lem}\label{Z-lem}
The identity map on $C_c(G\times G,B)\odot C^*(G)$ 
extends to an imprimitivity bimodule surjection
\begin{multline*}
\Gamma\:
{}_{(B\otimes C^*(G))\times G\times G\times N}
\bigl(Z_{G/N}^G((B\otimes C^*(G))\times_{\a\otimes\id}G)\bigr)_{
(B\otimes C^*(G))\times G\times G/N}
\\
\to
{}_{(B\times G\times G\times N)\otimes C^*(G)}
\bigl(Z_{G/N}^G(B\times_\a G)\otimes C^*(G)\bigr)_{
(B\times G\times G/N)\otimes C^*(G)}.
\end{multline*}
\end{lem}

\begin{proof}
The right-hand coefficient homomorphism
\[
\psi\:
(B\otimes C^*(G))\times_{\a\otimes\id} G
\times_{\hat{\a\otimes\id}|} G/N
\to(B\times_\a G\times_{\hat\a|} G/N)\otimes C^*(G)
\]
is the integrated form
$\psi = (k_B\otimes\id)\times(k_G\otimes 1)\times
(k_{C(G/N)}\otimes 1)$,
where $(k_B,k_G,k_{C(G/N)})$ are the canonical covariant
maps of $(B,G,C_0(G/N))$ into 
$M(B\times_\a G\times_{\hat\a}G/N)$.  
It is clear that $\psi$ maps generators to generators,
so is a surjection.

Similarly, the left-hand coefficient homomorphism
\[
\p\:
(B\otimes C^*(G))\times_{\a\otimes\id} G
\times_{\hat{\a\otimes\id}} G\times_{\hat{\hat{\a\otimes\id}}|} N
\to
(B\times_\a G\times_{\hat\a} G\times_{\hathata|} N)\otimes C^*(G)
\]
is the surjective integrated form
$\p = (\ell_B\otimes\id)\times (\ell_G\otimes 1)\times
(\ell_{C(G)}\otimes 1)\times (\ell_N\otimes 1)$,
where $(\ell_B, \ell_G, \ell_{C(G)}, \ell_N)$ are the 
canonical covariant maps of
$(B,G,C_0(G),N)$ into
$M(B\times_\a G\times_{\hat\a}G\times_{\hathata|}N)$.

For the bimodule map itself, a routine computation 
on elements of $C_c(G\times G,B)\odot C^*(G)$
shows that
\[
\<z\otimes c,w\otimes d\>_{(B\otimes C^*(G))\times G\times G/N}(s,tN)
=\<z,w\>_{B\times G\times G/N}(s,tN)\otimes c^*d,
\]
so we have
\[
\psi\bigl(
\<z\otimes c,w\otimes d\>_{(B\otimes C^*(G))\times G\times G/N}
\bigr)
=\<z,w\>_{B\times G\times G/N}\otimes c^*d.
\]
Therefore the identity map on $C_c(G\times G,B)\odot C^*(G)$ extends
to a bounded linear surjection
$\Gamma\:
Z_{G/N}^G((B\otimes C^*(G))\times_{\a\otimes\id}G)
\to
Z_{G/N}^G(B\times_\a G)\otimes C^*(G)$.

Another routine calculation (using \eqref{Z-gens}) shows that for
any generator $x$ from $B\otimes C^*(G)$, $G$, $C_0(G)$, or $N$, 
the left module actions on elementary tensors satisfy
\[
\Gamma(x\cdot(z\otimes c))
=\p(x)\cdot \Gamma(z\otimes c).
\]
Thus $\Gamma$ is a right-Hilbert bimodule 
homomorphism with coefficient maps $\p$ and~$\psi$.
\end{proof}

In particular, \lemref{Z-lem}
gives an embedding of $C_c(G\times G,B\otimes C^*(G))$ in 
$Z_{G/N}^G(B\times_\a G)\otimes C^*(G)$; by standard function space 
techniques (see for example \cite[Appendix C]{enchilada}),
we also get an 
embedding of $C_c(G\times G,M^\b(B\otimes C^*(G)))$ in 
$M(Z_{G/N}^G(B\times_\a G)\otimes C^*(G))$.
Thus (as in the proof of the next 
proposition), we can use the bimodule operations from 
$Z_{G/N}^G((B\otimes C^*(G))\times_{\a\otimes\id} G)$ 
in calculations with $C_c$-functions, 
although we regard the functions as elements of 
$Z_{G/N}^G(B\times_\a G)\otimes C^*(G)$.

\begin{prop}
Let $(B,G,\a)$ be an action, and let $N$ be a closed normal
subgroup of $G$.
The map $\d_Z\:C_c(G\times G,B)\to 
C_c(G\times G,M^\b(B\otimes C^*(G)))$
defined by
\[
\d_Z(z)(s,t)=z(s,t)\otimes t\inv s
\]
extends uniquely to an
$\infl\hat{\hathata|}-\hat\a\dec$
compatible coaction on 
$Z_{G/N}^G(B\times_\a G)$.
\end{prop}

\begin{proof}
Since $Z_{G/N}^G(B\times_\a G)$ is an imprimitivity bimodule, 
and since the coefficient maps are coactions,
it suffices to show that $\d_Z$
satisfies the coaction identity and
preserves the left module actions and the right inner products.
For the coaction identity, 
first note that for 
$y\otimes c\in C_c(G\times G,B)\odot C^*(G)$, we have
\begin{align*}
(\d_Z\otimes\id)(y\otimes c)
&\in C_c(G\times G,M^\b(B\otimes C^*(G))\odot C^*(G))
\\&\subset C_c(G\times G,M^\b(B\otimes C^*(G)\otimes C^*(G))),
\end{align*}
with
\begin{align*}
(\d_Z\otimes\id)(y\otimes c)(s,t)
&=(\d_Z(y)\otimes c)(s,t)
\\&=\d_Z(y)(s,t)\otimes c
\\&=y(s,t)\otimes t\inv s\otimes c
\\&=(\id\otimes\S)(y(s,t)\otimes c\otimes t\inv s)
\\&=(\id\otimes\S)\bigl((y\otimes c)(s,t)\otimes t\inv s\bigr).
\end{align*}
It follows from 
density and continuity for the inductive-limit topologies that
for $z\in C_c(G\times G,B)$,  we have
$(\d_Z\otimes\id)\circ\d_Z(z)
\in C_c(G\times G,M^\b(B\otimes C^*(G)\otimes C^*(G)))$,
with
\begin{align*}
(\d_Z\otimes\id)\circ\d_Z(z)(s,t)
&=(\id\otimes\S)\bigl(\d_Z(z)(s,t)\otimes t\inv s\bigr)
\\&=(\id\otimes\S)(z(s,t)\otimes t\inv s\otimes t\inv s)
\\&=z(s,t)\otimes t\inv s\otimes t\inv s
\\&=(\id\otimes\d_G)(z(s,t)\otimes t\inv s)
\\&=(\id\otimes\d_G)(\d_Z(z)(s,t))
\\&=(\id\otimes\d_G)\circ\d_Z(z)(s,t).
\end{align*}

For the left module actions, if $f\in C_c(N\times G\times G,B)
\subset B\times_a G\times_{\hat\a} G \times_{\hathata|}N$ and
$z\in C_c(G\times G,B)$ we have
\begin{align*}
&\Bigl(\infl\hat{\hathata|}(f)\cdot\d_Z(z)\Bigr)(s,t)
\\&\quad=\int_N\int_G\infl\hat{\hathata|}(f)(n,r,t)
(\a\otimes\id)_r\bigl(\d_Z(z)(r\inv s,r\inv tn)\bigr)
\D(n)^{1/2}\,dr\,dn
\\&\quad=\int_N\int_G(f(n,r,t)\otimes n)
(\a_r\otimes\id)\bigl(z(r\inv s,r\inv tn)\otimes n\inv t\inv s\bigr)
\D(n)^{1/2}\,dr\,dn
\\&\quad=\int_N\int_G(f(n,r,t)\otimes n)
\bigl(\a_r(z(r\inv s,r\inv tn))\otimes n\inv t\inv s\bigr)
\D(n)^{1/2}\,dr\,dn
\\&\quad=\int_N\int_G\bigl(f(n,r,t)\a_r(z(r\inv s,r\inv tn))
\otimes t\inv s\bigr)
\D(n)^{1/2}\,dr\,dn
\\&\quad=\int_N\int_G f(n,r,t)\a_r(z(r\inv s,r\inv tn))
\D(n)^{1/2}\,dr\,dn
\otimes t\inv s
\\&\quad=(f\cdot z)(s,t)\otimes t\inv s
=\d_Z(f\cdot z)(s,t).
\end{align*}

Finally, for the inner products, if $z,w\in C_c(G\times G,B)$ then
\begin{align*}
&\bigl\<\d_Z(z),\d_Z(w)\bigr\>_{
M((B\times G\times G/N)\otimes C^*(G)))}(s,tN)
\\&\quad=\int_N\int_G(\a\otimes\id)_{r\inv}
\bigl(\d_Z(z)(r,rtn)^*\d_Z(w)(rs,rtn)\bigr)\,dr\,dn
\\&\quad=\int_N\int_G(\a_{r\inv}\otimes\id)
\bigl((z(r,rtn)\otimes n\inv t\inv)^*
(w(rs,rtn)\otimes n\inv t\inv s)\bigr)\,dr\,dn
\\&\quad=\int_N\int_G(\a_{r\inv}\otimes\id)
\bigl(z(r,rtn)^*w(rs,rtn)\otimes s\bigr)\,dr\,dn
\\&\quad=\int_N\int_G \a_{r\inv}\bigl(z(r,rtn)^*w(rs,rtn)\bigr)\,dr\,dn
\otimes s
\\&\quad=\<z,w\>_{B\times G\times G/N}(s,tN)\otimes s
=\hat\a\dec\bigl(\<z,w\>_{B\times G\times G/N}\bigr)(s,tN).
\qedhere
\end{align*}
\end{proof}

\begin{thm}[Full Mansfield Imprimitivity]
\label{full mansfield}
For any maximal coaction $(A,G,\d)$
and any closed normal subgroup $N$ of $G$,
there exists an 
\[
A\times_\d G\times_{\hat\d|}N - A\times_{\d|}G/N
\]
imprimitivity bimodule $Y_{G/N}^G(A)$ 
and an $\infl\hat{\hat\d|}-\d\dec$ compatible
coaction $\d_Y$ of $G$ on $Y_{G/N}^G(A)$. 
\end{thm}

We refer to $Y_{G/N}^G(A)$ as the \emph{Mansfield bimodule}.

\begin{proof}
Let $(C,\epsilon) = (A\times_\d G\times_{\hat\d}G,\hathatd)$.
Then since $\delta$ is maximal, $K(A)$ is an equivariant $C-A$
imprimitivity bimodule, so the following diagram serves to
define $Y_{G/N}^G(A)$ and $\delta_Y$:
\begin{equation}\label{mansfield}
\xymatrix
@C+50pt
{
{(A\times_\d G\times_{\hat\d|} N, \infl\hat{\hat\d|})}
\ar[r]^-{(Y_{G/N}^G(A),\d_Y)}
&{(A\times_{\d|} G/N, \d\dec)}
\\
{(C\times_\e G\times_{\hat\e|} N, \infl\hat{\hat\e|})}
\ar[u]|{(K(A)\times G\times N,\infl\hat{\hat{\d_K}|})}
\ar[r]_-{(Z_{G/N}^G(C),\d_Z)}
&{(C\times_{\e|} G/N, \e\dec)}
\ar[u]|{(K(A)\times G/N,\d_K\dec)}.
}
\end{equation}
\end{proof}


\section{Naturality of the Mansfield Bimodule}
\label{natl-sec}

In this section we 
show that the Mansfield bimodule is a
natural equivalence between the
functors ``${}\times G\times N$'' and
``${}\times G/N$'', when those functors are restricted
to $\CC_m(G)$. 
The first step is to resolve two potential ambiguities
among our bimodules.

\begin{prop}\label{KZ-prop}
For any dual coaction $(A,G,\d)$, we have
\[
(K(A),\d_K)\cong(\zg(A),\d_Z).
\]
\end{prop}

\begin{proof}
Let $(B,G,\a)$ be an action such that
$(B\times_\a G,\hat\a) = (A,\d)$,
and define a linear map $\Psi\:C_c(G,B)\odot C_c(G)\to
C_c(G\times G,B)$ by
\[
\Psi(f\otimes\xi)(s,t)=f(s)\xi(t).
\]
To see that $\Psi$ extends to
an imprimitivity bimodule isomorphism
of $K(A)$ onto $Z_{G/G}^G(A)$,
it suffices to show that it preserves the right inner
products and the left module actions of the generators
of $B\times_\a G\times_{\hat\a}G \times_{\hathata}G$.
For the inner products, we have
\begin{align*}
&\<\Psi(f\otimes\xi),\Psi(g\otimes\eta)\>_{B\times G}(s)
\\&\quad=\iint\a_{t\inv}\bigl(\Psi(f\otimes\xi)(t,r)^*
\Psi(g\otimes\eta)(ts,r)\bigr)\,dr\,dt
\\&\quad=\iint\a_{t\inv}\bigl((f(t)\xi(r))^*g(ts)\eta(r)\bigr)\,dr\,dt
\\&\quad=\iint\a_{t\inv}\bigl(f(t)^*\bar{\xi(r)}
g(ts)\eta(r)\bigr)\,dr\,dt
\\&\quad=\iint\a_{t\inv}\bigl(f(t)^*g(ts)\bigr)
\bar{\xi(r)}\eta(r)\,dr\,dt
\\&\quad=\int\a_{t\inv}\bigl(f(t)^*g(ts)\bigr)\,dt
\int\bar{\xi(r)}\eta(r)\,dr
\\&\quad=(f^*g)(s)\<\xi,\eta\>
=\<f,g\>_{B\times G}(s)\<\xi,\eta\>
\\&\quad=\bigl(\<f,g\>_{B\times G}\<\xi,\eta\>\bigr)(s)
=\<f\otimes\xi,g\otimes\eta\>_{B\times G}(s).
\end{align*}

Recall that the left action on $K(A)$ is implemented by the 
canonical surjection
$\Phi = (\id\otimes\lambda)\circ\hat\a \times 
(1\otimes M) \times (1\otimes\rho)$; the 
left action of the generators on $Z_{G/G}^G(A)$ 
is given by \eqref{Z-gens}.  
Thus, if $(i_B,i_G)$ denotes the canonical 
covariant homomorphism of $(B,G)$ into $M(B\times_\alpha G)$,
for $b\in B$ we have
\begin{align*}
&\Psi\bigl(
b\cdot(f\otimes\xi)
\bigr)(s,t)
\\&\quad=\Psi\bigl(
(\id\otimes\l)\circ\hat\a(i_B(b))
(f\otimes\xi)
\bigr)(s,t)
\\&\quad=\Psi\bigl(
(i_B(b)\otimes 1)
(f\otimes\xi)
\bigr)(s,t)
\\&\quad=\Psi(i_B(b)f\otimes\xi)(s,t)
=\bigl(i_B(b)f\bigr)(s)\xi(t)
\\&\quad=bf(s)\xi(t)
=b\Psi(f\otimes\xi)(s,t)\\
&\quad=\bigl(b\cdot\Psi(f\otimes\xi)\bigr)(s,t),
\end{align*}
and for $r\in G$ we have
\begin{align*}
&\Psi\bigl(
r\cdot(f\otimes\xi)
\bigr)(s,t)
\\&\quad=\Psi\bigl(
(\id\otimes\l)\circ\hat\a(i_G(r))
(f\otimes\xi)
\bigr)(s,t)
\\&\quad=\Psi\bigl(
(i_G(r)\otimes\l_r)
(f\otimes\xi)
\bigr)(s,t)
\\&\quad=\Psi(i_G(r)f\otimes\l_r\xi)(s,t)
=\bigl(i_G(r)f\bigr)(s)(\l_r\xi)(t)
\\&\quad=\a_r\bigl(f(r\inv s)\bigr)\xi(r\inv t)
=\a_r\bigl(\Psi(f\otimes\xi)(r\inv s,r\inv t)\bigr)
\\&\quad=\bigl(r\cdot\Psi(f\otimes\xi)\bigr)(s,t).
\end{align*}
Similarly,
for $h\in C_c(G)\subset C_0(G)$ we have
\begin{align*}
\Psi\bigl(
h\cdot(f\otimes\xi)
\bigr)(s,t)
&=\Psi\bigl(
(1\otimes M_h)(f\otimes\xi)
\bigr)(s,t)
\\&=\Psi(f\otimes M_h\xi)(s,t)
\\&=f(s)(M_h\xi)(t)
=f(s)h(t)\xi(t)
\\&=h(t)f(s)\xi(t)
=h(t)\Psi(f\otimes\xi)(s,t)
\\&=\bigl(h\cdot\Psi(f\otimes\xi)\bigr)(s,t),
\end{align*}
and for $u\in G$ 
(meaning the ``outermost'' $G$ in 
$B\times_\a G\times_{\hat\a}G\times_{\hathata}G$)
we have
\begin{align*}
\Psi\bigl(
u\cdot(f\otimes\xi)
\bigr)(s,t)
&=\Psi\bigl(
(1\otimes\rho_u)
(f\otimes\xi)
\bigr)(s,t)
\\&=\Psi(f\otimes\rho_u\xi)(s,t)
=f(s)(\rho_u\xi)(t)
\\&=f(s)\xi(tu)\D(u)^{1/2}
=\Psi(f\otimes\xi)(s,tu)\D(u)^{1/2}
\\&=\bigl(u\cdot\Psi(f\otimes\xi)\bigr)(s,t).
\end{align*}

Finally, we show $\Psi$ is $\d_K-\d_Z$ equivariant.
Let $f\in C_c(G,B)\subset B\times_\a G$ and $\xi\in C_c(G)\subset L^2$.
Then
\begin{align}
(\Psi\otimes\id)\circ\d_K(f\otimes\xi)(s,t)
&=(\Psi\otimes\id)
\bigl(W(\hat\a\otimes_*\id)(f\otimes\xi)\bigr)(s,t)
\notag
\\&=(\Psi\otimes\id)
\bigl(w_G^*\cdot(\id\otimes\S)(\hat\a(f)\otimes\xi)\bigr)(s,t)
\notag
\\&=\bigl(w_G^*\cdot(\Psi\otimes\id)\circ
(\id\otimes\S)(\hat\a(f)\otimes\xi)\bigr)(s,t)
\notag
\\&=\bigl(1\otimes w_G^*(t)\bigr)\bigl((\Psi\otimes\id)\circ
(\id\otimes\S)(\hat\a(f)\otimes\xi)\bigr)(s,t)
\notag
\\&=(1\otimes t\inv)\hat\a(f)(s)\xi(t)
\tag{*}
\\&=(1\otimes t\inv)\bigl(f(s)\otimes s\bigr)\xi(t)
\notag
\\&=f(s)\xi(t)\otimes t\inv s
\notag
\\&=\Psi(f\otimes\xi)(s,t)\otimes t\inv s
\notag
\\&=\d_Z\circ\Psi(f\otimes \xi)(s,t),
\notag
\end{align}
where the equality (*) is verified by temporarily replacing $\hat\a(f)$
by an elementary tensor $g\otimes c\in C_c(G,B)\odot C^*(G)$,
computing that
\begin{align*}
(\Psi\otimes\id)\circ(\id\otimes\S)(g\otimes c\otimes \xi)(s,t)
&=(\Psi\otimes\id)(g\otimes \xi\otimes c)(s,t)
\\&=\bigl(\Psi(g\otimes \xi)\otimes c\bigr)(s,t)
\\&=\Psi(g\otimes \xi)(s,t)\otimes c
\\&=g(s)\xi(t)\otimes c
\\&=(g(s)\otimes c)\xi(t)
\\&=(g\otimes c)(s)\xi(t),
\end{align*}
and then appealing to continuity and density for the inductive-limit
topologies.
\end{proof}

\begin{cor}\label{YK-cor}
For any maximal coaction $(A,G,\d)$, 
we have
\[
(\yg(A),\d_Y)\cong(K(A),\d_K).
\]
\end{cor}

\begin{proof}
Applying naturality of Katayama (\thmref{natural katayama})
to the Katayama bimodule itself gives
\[
K(A)\times_{\d_K}G\times_{\hat{\d_K}}G
\otimes_{A\times G\times G} K(A)
\cong
K(A\times_{\d}G\times_{\hat{\d}}G)
\otimes_{A\times G\times G} K(A),
\]
equivariantly for the appropriate coactions.  
Since $\delta$ is maximal, $K(A)$ is an imprimitivity bimodule, 
so it follows that
$K(A)\times_{\d_K}G\times_{\hat{\d_K}}G$ is equivariantly
isomorphic to 
$K(A\times_{\d}G\times_{\hat{\d}}G)$.
Combining this with \propref{KZ-prop} (applied to $\hathatd$)
gives
\[
K(A)\times_{\d_K}G\times_{\hat{\d_K}}G
\cong
Z_{G/G}^G(A\times_{\d}G\times_{\hat{\d}}G),
\]
equivariantly for the appropriate coactions; the proposition
now follows immediately from this and the definition 
of $Y_{G/G}^G(A)$ (Diagram~\eqref{mansfield}). 
\end{proof}

We now establish a curious duality between Mansfield and Green
imprimitivity. 

\begin{prop}[Factorization of $\zn$]
\label{Z-fact}
For any dual coaction $(A,G,\d)$ and any closed normal
subgroup $N$ of $G$, the diagram
\[
\xymatrix{
{(A\times_\d G\times_{\hat\d|} N,\infl\hat{\hat\d|})}
\ar[rr]^-{(\zn(A),\d_Z)}
&&{(A\times_{\d|} G/N,\d\dec)}
\\
&{(A\times_\d G\times_{\hat\d} G\times_{\hathatd|} G/N,
\hathatd{}\dec)}
\ar[ul]|{(\xn(A\times G),\d_X)}
\ar[ur]|{(K(A)\times G/N,\d_K\dec)}
}
\]
commutes.
\end{prop}

\begin{proof}
Let $(B,G,\a)$ be an action such that 
$(B\times_\a G,\hat\a) = (A,\d)$.  Then
the desired diagram is the upper left triangle of the diagram
\[
\xymatrix
@C+60pt
{
{(A\times_\d G\times_{\hat\d} G\times_{\hathatd|}G/N, 
 \hathatd{}\dec)}
\ar[r]^-{(\xn(A\times G),\d_X)}
\ar[d]|{(K(A)\times G/N,\d_K\dec)}
&{\bigl(B\times_\a G\times_{\hat\a} G\times_{\hathata|} N,
 \infl\hat{\hathata|}\bigr)}
\ar[dl]|{(\zn(B\times G),\d_Z)}
\ar[d]|{\bigl(\xg(B)\times N,\infl\hat{\a^X|}\bigr)}
\\
{(B\times_\a G\times_{\hat\a|} G/N,\hat\a\dec)}
\ar[r]_{(\xn(B),\d_X)}
&{(B\times_{\a|} N,\infl\hat{\a|})},
}
\]
where $\a^X$ is the 
$\hathata-\a$ compatible action on $\xg(B)$ 
from \cite[Theorem 1]{ech:twisted}.
Without the equivariance,
the lower right triangle commutes by \cite[Theorem 3.1]{ekr}.
Equivariance is proven in \cite[Theorem 5.1]{enchilada}, albeit in the
context of reduced crossed products; the argument given there carries
over with no significant change, since it only uses $C_c$-functions and
pre-right-Hilbert bimodules.

Next, note that the lower triangle in the special case
$N=G$ gives us 
\[
(\zg(B\times_\a G),\d_X)
\cong (\xg(B)\times_{\alpha^X} G,\hat{\a^X}),
\]
so by equivariance and \propref{KZ-prop} we get
\[
\bigl(K(A)\times_{\d_K|} G/N,\d_K\dec\bigr)
\cong 
\bigl(\xg(B)\times_{\a^X} G\times_{\hat{\a^X}|} G/N,
 \hat{\a^X}{}\dec\bigr).
\]
Using this to replace the left-hand vertical arrow in
the diagram above
reduces the outer rectangle to naturality of Green
(\thmref{natural green}),
applied to the $A\times_\d G - B$ imprimitivity bimodule
$\xg(B)$;
this finishes the proof.
\end{proof}

\begin{cor}[Factorization of $\yn$]
\label{Y-fact}
For any maximal coaction $(A,G,\d)$ and any 
closed normal subgroup $N$ of $G$, the diagram
\[
\xymatrix{
{(A\times_\d G\times_{\hat\d|} N,\infl\hat{\hat\d|})}
\ar[rr]^-{(\yn(A),\d_Y)}
&&{(A\times_{\d|} G/N,\d\dec)}
\\
&{(A\times_\d G\times_{\hat\d} G\times_{\hathatd|} G/N,
\hathatd{}\dec)}
\ar[ul]|{(\xn(A\times G),\d_X)}
\ar[ur]|{(K(A)\times G/N,\d_K\dec)}
}
\]
commutes.
\end{cor}

\begin{proof}
To simplify the notation, let 
$(C,G,\e) = (A\times_\d G\times_{\hat\d} G,G,\hathatd)$. 
Then the desired diagram is the top triangle of the diagram
\[
\xymatrix{
{A\times_\d G\times_{\hat\d|} N}
\ar[rr]^-{\yn(A)}
&&{A\times_{\d|} G/N}
\\
&{A\times_\d G\times_{\hat\d} G\times_{\hathatd|} G/N}
\ar[ul]|{\xn(A\times G)}
\ar[ur]|{K(A)\times G/N}
\\
&{C\times_\e G\times_{\hat\e}G\times_{\hathate|} G/N}
\ar[dl]|{\xn(C\times G)}
\ar[dr]|{K(C)\times G/N}
\ar[u]|{K(A)\times G\times G\times G/N}
\\
{C\times_\e G\times_{\hat\e|} N}
\ar[rr]_-{\zn(C)}
\ar[uuu]|{K(A)\times G\times N}
&&{C\times_{\e|} G/N},
\ar[uuu]|{K(A)\times G/N}
}
\]
which is to be interpreted as an equivariant diagram,
although we have omitted the appropriate coactions for clarity.
The bottom triangle commutes by \propref{Z-fact}.  The outer
rectangle commutes by definition of $Y_{G/N}^G(A)$
(Diagram~\eqref{mansfield}).  
The left quadrilateral
commutes by naturality of Green 
(\thmref{natural green}),
applied to the $C\times_\e G - A\times_\d G$ imprimitivity
bimodule $K(A)\times_{\d_K}G$.  
Finally, the right quadrilateral
commutes by naturality of Katayama
(\thmref{natural katayama}),
applied to $K(A)$ itself,
together with functoriality of ``${}\times G/N$''.
\end{proof}

We can use the above result to clear up another apparent ambiguity,
which arises for dual coactions.

\begin{prop}
For any dual coaction $(A,G,\d)$ and any closed normal
subgroup $N$ of $G$, we have
\[
(\yn(A),\d_Y)\cong (\zn(A),\d_Z).
\]
\end{prop}

\begin{proof}
This is immediate from \propref{Z-fact}
and \corref{Y-fact}, since both bimodules
factor as
\[
\widetilde{\xn}(A\times G)
\otimes_{A\times_\d G\times_{\hat\d}G\times G/N}
(K(A)\times G/N),
\]
equivariantly for the appropriate coactions.
\end{proof}

We are now ready to prove that the Mansfield bimodule is natural.
We have factored Mansfield as a composition of Green and Katayama,
each of which is natural, and this is what makes it work.

\begin{thm}[Naturality of Mansfield]
\label{natural mansfield}
Let ${}_{(A,\d)}(X,\z)_{(B,\e)}$ be a right-Hilbert bimodule coaction
of a locally compact group $G$, and let $N$ be a closed normal
subgroup of $G$.
If $\d$ and $\e$ are maximal, then the diagram
\[
\xymatrix
@C+30pt
{
{(A\times_\d G\times_{\hat\d|}N,\infl\hat{\hat\d|})}
\ar[r]^-{(\yn(A),\d_Y)}
\ar[d]_{(X\times G\times N,\infl\hat{\hat\z|})}
&{(A\times_{\d|}G/N,\d\dec)}
\ar[d]^{(X\times G/N,\z\dec)}
\\
{(B\times_\e G\times_{\hat\e|}N,\infl\hat{\hat\e|})}
\ar[r]_-{(\yn(B),\d_Y)}
&{(B\times_{\e|}G/N,\e\dec)}
}
\]
commutes.
\end{thm}

\begin{proof}
The desired diagram is the outer rectangle of the diagram
\[
\xymatrix{
{A\times_\d G\times_{\hat\d|}N}
\ar[rr]^-{\yn(A)}
\ar[ddd]_{X\times G\times N}
&&{A\times_{\d|}G/N}
\ar[ddd]^{X\times G/N}
\\
&{A\times_\d G\times_{\hat\d} G\times_{\hathatd|} G/N}
\ar[ul]|{\xn(A\times G)}
\ar[ur]|{K(A)\times G/N}
\ar[d]|{X\times G\times G\times G/N}
\\
&{B\times_\e G\times_{\hat\e} G\times_{\hathate|} G/N}
\ar[dl]|{\xn(B\times G)}
\ar[dr]|{K(B)\times G/N}
\\
{B\times_\e G\times_{\hat\e|}N}
\ar[rr]_-{\yn(B)}
&&{B\times_{\e|}G/N},
}
\]
which is again to be interpreted equivariantly  for the 
appropriate coactions.
The top and bottom triangles commute by factorization of $\yn$ 
(\corref{Y-fact}),
the left quadrilateral commutes by naturality of Green
(\thmref{natural green}), 
and the right
quadrilateral commutes by naturality of Katayama 
(\thmref{natural katayama})
and functoriality of ``${}\times G/N$''.
\end{proof}


\section{Preservation of maximality}
\label{sec:preserve}

In this section, partly as an application of the results of the
preceding two, we show that the functors involved
in the Full Mansfield Imprimitivity Theorem 
(\thmref{full mansfield}) can be regarded as having 
the category $\CC_m(G)$
of maximal coactions of $G$ not only for their domain, but also
their codomain.  This involves verifying that maximality of coactions
is preserved under various operations; in doing so we freely
use the facts that dual coactions are maximal,
and that maximality is preserved by Morita-Rieffel equivalence
(\cite[Propositions~3.4 and 3.5]{ekq}).
We end the section with another
preservation-of-maximality result which is not needed for the purpose
of investigating the codomains of our functors, but which is of obvious
interest in this circle of ideas.

First recall that if $(A,G,\a,N,\tau)$ is a twisted action in the sense
of Green \cite{gre:local}, then there is a natural ``dual coaction''
$\tilde\a$ of $G/N$ on the twisted crossed product 
$A\times_{\a,N}G$.

\begin{prop}\label{dual-max}
Let $(A,G,\a,N,\tau)$ be a twisted action. Then the dual coaction
$(A\times_{\a,N}G,G/N,\tilde\a)$ is maximal.
\end{prop}

\begin{proof}
By Echterhoff's version of the Stabilization Trick \cite[Theorem
1]{ech:twisted}, the twisted action of $(G,N)$ is Morita-Rieffel
equivalent to an inflated action $(B,G,\infl\b)$, where $\b$ is an
action of $G/N$ on $B$ (and this inflated action is trivially twisted
over~$N$).  Thus the dual coactions $\hat\a$ and $\hat{\infl\b}$ are
Morita-Rieffel equivalent, hence so are the 
corresponding coactions $\tilde\a$
and $\tilde{\hat{\infl\b}}$ on the twisted crossed products. But the
twisted crossed product 
$B\times_{\infl\beta,N}G$ 
is isomorphic to the ordinary crossed
product $B\times_\b G/N$, and $\tilde{\hat{\infl\b}}$
corresponds to the dual coaction $\hat\b$ under this
isomorphism.  Consequently, 
$\tilde\a$ is Morita-Rieffel equivalent to $\hat\b$;
since $\hat\b$ is maximal, so is $\tilde\a$.
\end{proof}

\begin{cor}
Let $(A,G,\d)$ be a maximal coaction, and let $N$ be a closed normal
subgroup of $G$. Then the restricted coaction 
$(A,G/N,\d|)$ is also maximal.
\end{cor}

\begin{proof}
By maximality, $\d$ is Morita-Rieffel equivalent 
(via the Katayama bimodule of \propref{kat-coact}) 
to a dual coaction
$(B\times_\a G,\hat\a)$; hence the restricted coactions $\d|$
and $\hat\a|$ are also Morita-Rieffel equivalent.  By Green's
Decomposition Theorem \cite[Proposition 1]{gre:local}, the ordinary
crossed product $B\times_\a G$ is isomorphic to the twisted crossed
product $B\times_{\a|}N\times_{\a\dec,N}G$, and it is easy to see that,
under this isomorphism,
$\hat\a|$ corresponds to $\tilde{\a\dec}$,
which is maximal by \propref{dual-max}.
\end{proof}

\begin{prop}\label{infl-max}
Let $(A,H,\d)$ be a maximal coaction, where $H$ is a closed subgroup of
a locally compact group $G$. 
Then the inflated coaction $(A,G,\infl\d)$ is also maximal.
\end{prop}

\begin{proof}
By maximality, $\d$ is Morita-Rieffel equivalent 
(via the Katayama bimodule of \propref{kat-coact}) 
to a dual coaction
$(B\times_\a H,\hat\a)$; hence the inflated coactions $\infl\d$ and
$\infl\hat\a$ of $G$ are also Morita-Rieffel equivalent.
By Green's Imprimitivity Theorem for induced algebras
\cite[Theorem~17]{gre:local},
$B\times_\a H$
is Morita-Rieffel equivalent to
the crossed product
$\ind_H^GA\times_{\ind\a}G$ by the induced action,
and it is easy to see that
the inflated coaction $\infl\hat\a$
corresponds to
the dual coaction $\hat{\ind\a}$ under this isomorphism.
\end{proof}

\begin{prop}
Let $(A,G,\d)$ be a maximal coaction, and let $N$ be a closed normal
subgroup of $G$. Then the decomposition coaction 
$(A\times_{\d|}G/N, G, \d\dec)$ 
is also maximal.
\end{prop}

\begin{proof}
By maximality, $\d\dec$ is Morita-Rieffel equivalent 
(via the Mansfield bimodule of \thmref{full mansfield}) 
to $\infl\hat{\hat\d|}$, which is maximal by \propref{infl-max}.
\end{proof}

Let $(B,G,\e)$ be a coaction. 
A (closed two-sided) ideal $J$ of $B$ 
is {invariant} under $\e$ in the sense of 
\cite{nil:full}  if 
\[
\bar{\e(J)(1\otimes C^*(G))}=J\otimes C^*(G);
\]
Nilsen has shown \cite[Proposition 2.1]{nil:full} that 
when this is the case, the restriction $\e|_J$ is a coaction 
of $G$ on $J$. 

\begin{prop}\label{invt-max}
Let $(B,G,\e)$ is a maximal coaction,
and let $J$ be an ideal of $B$
which is invariant under $\e$ in the above sense. 
Then the restriction $(J,G,\e|_J)$ is also maximal.
\end{prop}

\begin{proof}
Let $(\pi,\mu)$ be the canonical covariant
homomorphism of $(B,C_0(G))$ into $M(B\times_\e G)$; then
by \cite[Proposition 2.1]{nil:full},
$\pi|_J\times\mu$ is an isomorphism of $J\times_{\e|_J} G$ 
onto an ideal of $B\times_\e G$.%
\footnote{There is a subtlety concerning how to make an
``official'' covariant homomorphism out of $(\pi|_J,\mu)$, but Nilsen
shows that this causes no problem.} 
Identifying $J\times_{\e|_J} G$ with its image in 
$B\times_\e G$,
it is easy to see that this ideal is invariant under the dual
action $\hat\e$. Thus by \cite[Proposition 12]{gre:local}, the double
crossed product $J\times_{\e|_J} G\times_{\hat{\e|_J}} G$ 
can be identified with an
ideal of $B\times_\e G\times_{\hat\e} G$. 
But subject to these identifications, 
the canonical surjection 
$\Phi_J\:J\times_{\e|_J} G\times_{\hat{\e|_J}} G\to
J\otimes\K$ is just the restriction of the isomorphism 
$\Phi_B\:B\times_\e G\times_{\hat\e} G\to B\otimes\K$,
and hence is 1-1.
\end{proof}

Associated to a right-Hilbert bimodule coaction
${}_{(A,\d)}(X,\z)_{(B,\e)}$, there can be as many as three other
$C^*$-coactions: first, 
there is a coaction $\mu$ on the imprimitivity algebra $\K_B(X)$ such
that $\zeta$ is $\mu-\e$ compatible, 
by \cite[Proposition 2.30]{enchilada};
second, there is a coaction $\nu$
on the linking algebra $L(X)$ which compresses on the corners to $\mu$,
$\zeta$, and $\e$,
by \cite[Lemma~2.22]{enchilada};
third, 
the closed span $B_X=\bar{\<X,X\>_B}$ is an $\e$-invariant
ideal of $B$, and the restriction
of $\e$ to $B_X$ 
is a coaction $\vartheta$ on $B_X$
such that $\zeta$ is $\mu-\vartheta$ compatible,
by \cite[Lemma 2.32]{enchilada}.

\begin{cor}
With notation as above, if
the coaction $(B,G,\e)$ is maximal, then so are
$(\K_B(X),G,\mu)$, 
$(L(X),G,\nu)$, and
$(B_X,G,\vartheta)$.
\end{cor}

\begin{proof}
Since $\mu$, $\nu$, and $\vartheta$ are all Morita-Rieffel
equivalent, it suffices to show that $\vartheta$ is maximal; 
but this is immediate from \propref{invt-max}. 
\end{proof}


\providecommand{\bysame}{\leavevmode\hbox to3em{\hrulefill}\thinspace}
\providecommand{\MR}{\relax\ifhmode\unskip\space\fi MR }
\providecommand{\MRhref}[2]{%
  \href{http://www.ams.org/mathscinet-getitem?mr=#1}{#2}
}
\providecommand{\href}[2]{#2}

\end{document}